\pdfoutput=1
\RequirePackage{ifpdf}
\ifpdf % We are running pdfTeX in pdf mode
\documentclass[pdftex]{sigma}
\else
\documentclass{sigma}
\fi

\begin{document}

\allowdisplaybreaks

\newcommand{\arXivNumber}{1910.12864}

\renewcommand{\thefootnote}{}

\renewcommand{\PaperNumber}{024}

\FirstPageHeading

\ShortArticleName{Horospherical Cauchy Transform on Some Pseudo-Hyperbolic Spaces}

\ArticleName{Horospherical Cauchy Transform\\ on Some Pseudo-Hyperbolic Spaces\footnote{This paper is a~contribution to the Special Issue on Algebra, Topology, and Dynamics in Interaction in honor of Dmitry Fuchs. The full collection is available at \href{https://www.emis.de/journals/SIGMA/Fuchs.html}{https://www.emis.de/journals/SIGMA/Fuchs.html}}}

\Author{Simon GINDIKIN}

\AuthorNameForHeading{S.~Gindikin}

\Address{Department of Mathematics, Hill Center, Rutgers University,\\
110 Frelinghysen Road, Piscataway, NJ 08854, USA}
\Email{\href{mailto:sgindikin@gmail.com}{sgindikin@gmail.com}}

\ArticleDates{Received October 28, 2019, in final form March 29, 2020; Published online April 07, 2020}

\Abstract{We consider the horospherical transform and its inversion in 3 examples of hyperboloids. We want to illustrate via these examples the fact that the horospherical inversion formulas can be directly extracted from the classical Radon inversion formula. In a~more broad context, this possibility reflects the fact that the harmonic analysis on symmetric spaces (Riemannian as well as pseudo-Riemannian ones) is equivalent (homologous), up to the Abelian Fourier transform, to the similar problem in the flat model. On the technical level it is important that we work not with the usual horospherical transform, but with its Cauchy modification.}


\Classification{32A45; 33C55; 43A75; 44A12}

\begin{flushright}
\begin{minipage}{70mm}
\it To Mitya Fuchs, my dear friend of more\\ than 60 years, on his 80th birthday.
\end{minipage}
\end{flushright}

\renewcommand{\thefootnote}{\arabic{footnote}}
\setcounter{footnote}{0}

\section{Introduction}

We realize pseudo-hyperbolic geometries on the $(n-1)$-dimensional hyperboloids $X_{p,q} \subset {\mathbb R}^n$, $n=p+q,$ defined by the equations
\begin{gather*}
\square(x)=\square _{p,q} (x)=(x_1)^2+\cdots+(x_p)^2- (x_{p+1})^2-\cdots -(x_{p+q})^2=1.
\end{gather*}
We do not exclude the case $q=0$, where we have the sphere. $X_{p,q}$ carries a transitive action of the group ${\rm SO}(p,q)$, and this action preserves the pseudo-hyperbolic metric.

We will consider the following 3 examples. In the 1st example $p=1$, $q=n-1$, we have the classical hyperbolic geometry, and this is an example of a~Riemannian non-compact symmetric space. The 2nd example is that of the sphere $S^{n-1}=X_{n,0}$. This is an example of a~compact Riemannian symmetric space. Our final example is that of $X_{2,n-2}$, which models the pseudo-hyperbolic geometry; it is a pseudo-Riemannian symmetric space that includes the group ${\rm SL}(2;{\mathbb R})$.

In all of these cases the inversion of the horospherical transform is the result of the same construction: the application of a universal fundamental closed differential form on the set of hyperbolic sections of hyperboloids. We deform the cycle of horospheres into a special cycle of hyperplane sections (we call it a geodesic cycle), on which the inverse transform is equivalent to the one given by the classical Radon inversion formula.

Let us discuss the conceptual picture which we illustrate by the simplest examples in this paper. There are 2 parallel languages in harmonic analysis on symmetric spaces. One is the original spectral approach that uses the decomposition into irreducible representations. This can be interpreted as the spherical Fourier transform. The other one is the horospherical transform of Gelfand -- the analog of the affine Radon transform in which affine hyperplanes are replaced by horospheres in the symmetric
spaces. These 2 transforms are connected by the (Abelian) Mellin transform. The translation from one language to another can help with some problems. Here is the most important example: the inverse spherical Fourier transform (Plancherel formula) can be deduced from the inverse horospherical formula.

It may not be easy to use this simple connection between 2 languages in order to obtain explicit statements and formulas as there are facts that are specific to one of the languages.
The focus of this paper is one such fact in the example of hyperboloids: the inversion of the horospherical transform is equivalent (homologous!) to a similar problem for the flat models, which is a simple modification of the Radon inversion formula.

Specifically, we consider for functions with support on a hyperboloid their  integrals over hyperplane sections. The problem of the reconstruction of functions through these integrals is overdetermined. Our principal construction is that of a universal closed differential form on the space of hyperplane sections. We call it the fundamental form for hyperboloids. Different inversion formulas for the Radon transform on the hyperboloids can be produced by the integration of this closed form over different cycles in the space of sections.

A special role is played by the cycle which we call geodesic. It is the cycle of hyperplane sections passing through the origin. Using the central projection we can interpret the restriction to the geodesic cycle as a version of the projective Radon transform~\cite{3}, and the Radon inversion formula gives the inversion of this transform. It plays the role of the flat model for our problem. In the case of the sphere, it is the Minkowski--Funk transform.

The inversion formula on the geodesic cycle coincides with the restriction of the fundamental closed form. The same will be true for homologous cycles. In particular, we construct a contraction of the cycle of horospheres to the geodesic cycle. In the case of a compact or a~pseudo-Riemannian symmetric space we need one more essential ingredient that we obtain by defining the complex horospherical transform on real hyperboloids. As we illustrate by the 3rd example, we need the complex horospheres if there are discrete series of representations. The construction of the horospherical transform is more complicated and takes values in $\bar \partial$-cohomology. We will consider it in another paper.

\section{Preliminary constructions}
Let us start with some definitions and facts for an arbitrary $X=X_{p,q}$. We will assume here that all functions satisfy $f\in C_0 ^\infty(X)$. Sometimes we interpret them as functions on ${\mathbb R}^n$ with support on the hyperboloid~$X$. We will need
some notations for differential forms. Let us denote by
$[a_1,\dots,a_n]$ the determinant of the matrix with columns
$a_1,\dots,a_n$, some of the entries can be 1-forms. We expand such
determinants from left to right and use the exterior product for
the multiplication of 1-forms. Such a determinant with identical
columns can differ from zero: $[{\rm d}x,\dots,{\rm d}x]=n! {\rm d}x_1\wedge\cdots
\wedge {\rm d}x_n$. We will write $a^{\{k\}}$ if a~column~$a$ is repeated~$k$ times.

Recall that by definition the interior product of forms, $\varphi \rfloor\psi$,
is a form $\alpha$ such that $ \varphi \wedge\alpha=\psi$. Its
restriction to the submanifold where $\varphi =0$ is uniquely
defined. If $\varphi={\rm d} f$ where $f$ is a~function, then ${\rm d} f\rfloor
\psi$ is, up to a constant factor, the residue of $\psi/f$ at $\{f=0\}$.

Let ${\mathbb R}^n_\xi$ and ${\mathbb R}^n_x$ be two copies of ${\mathbb R}^n$, and consider the pairing
\[
{\mathbb R}^n_\xi\times{\mathbb R}^n_x\longrightarrow {\mathbb R},\qquad (\xi,x)\mapsto \langle \xi,x\rangle,
\]
where the bilinear form $\langle \xi,x\rangle $ is the one that corresponds to the quadratic form $\square _{p,q}$. Let us define the Radon--Cauchy transform \cite{1,5} of $f$ to be
\[ \hat f(\xi,p)=\int _X \frac {f(x)}{\langle \xi,x \rangle -p-{\rm i}0} \omega_{p,q}(x,{\rm d}x),\]
where
\[ \omega_{p,q}(x,{\rm d}x)=2{\rm d} (\square_{p,q} (x))\rfloor
\big[{\rm d}x^{\{n\}}\big]=(n-1)!\sum _{1\leq j\leq n} (-1)^{(j-1)}\delta_jx_j
\bigwedge_{i\neq j} {\rm d}x_i ;\]
 here $\delta_j$ are the coefficients of $(x_j)^2$ in $\square _{p,q}(x)$. It is invariant on $X$. Also
 \[ \lim_ {\varepsilon \rightarrow 0} \frac 1 {p-{\rm i}\varepsilon} =(p-{\rm i} 0)^{-1}=p^{-1}+{\rm i}\pi\delta(p).
 \]

This definition is obtained from the usual definition of the Radon transform by replacing~$\delta(t)$, $t=\langle \xi,x \rangle -p$, with the distribution $(t-{\rm i}0)^{-1}$. As a result, the usual Radon inversion formula, originally valid only for~$n$ odd, gets extended to any~$n$. The Radon and Cauchy--Radon transforms are easily expressed through one another. We have
\[ \hat f(\lambda \xi ,\lambda p)=\lambda^{-1} \hat f(\xi, p), \qquad \lambda >0.\]
Let $\hat f(\xi)=\hat f(\xi, 1)$.

So the Cauchy--Radon transform on $X_{p,q}$ is a special case of the Cauchy--Radon transform in ${\mathbb R}^n$. We investigate specific properties of $\hat f$ if $f\in C_0^\infty(X)$ (i.e., has the support in~$X$). Let us remark that $\dim X=n-1$, but $\hat f$ is the result of integration on an $n$-parametric family of sections of $X$ by the hyperplanes in ${\mathbb R}^n$. Specifically, this can be expressed as the fact that $\hat f$ is a solution of the ultrahyperbolic equation
\[ \left\{\square _{p,q}\left(\frac {\partial}{\partial \xi}\right)-\frac {\partial^2}{\partial p^2}\right\} \hat f(\xi,p)=0.
\]

The problem of reconstructing $f$ on $X$ from $\hat f$ is overdetermined and equivalent to boundary problems for these ultrahyperbolic equations. Therefore it is natural to reconstruct $f$ on $X$ not from $\hat f $ for all $\xi$, $p$, but only from some $(n-1)$-parametric subfamilies. Of course, we need to remember that $\hat f$ is homogeneous.

Finally, recall the following technical lemma
\begin{lemma}
\[  {\rm d} \big[a(\xi), \xi, {\rm d}\xi^{\{n-2\}}\big]=-\frac 1
{n-2}\left(\sum_{1\leq j\leq n} \frac {\partial a_j(\xi)}{\partial \xi_j}\right)\big[\xi,{\rm d}\xi^{\{n-1\}}\big].
\]
\end{lemma}
The proof is a direct computation; it is simplified by the observation that it is sufficient to consider the case where the column~$a(\xi)$ has at most one non-zero element.

\section{Fundamental form for hyperboloids} As we mentioned above, functions $f$ on hyperboloids depend on $n-1$ variables, but their Radon transform $\hat f$ depends on~$n$ variables. Therefore the inversion formula cannot be unique. We will relate various such formulas
 with certain $(n-1)$-dimensional subfamilies of hyperplane sections.

 Our basic tool is a remarkable fundamental differential form on $X\times {\mathbb R}^n_\xi$:
\begin{gather*} \kappa_x[f]=\frac {f(u)}{(\langle \xi , (u-x)\rangle -{\rm i} \varepsilon)^{n-1}}
\omega_{p,q}(u,{\rm d}u)\wedge \big[x+u,\xi,{\rm d}\xi ^{\{n-2\}}\big],\qquad u\in X
 ,\quad \xi \in {\mathbb R}^n\backslash \{0\}.
\end{gather*}
 Here $f$ is a fixed function on the hyperboloid $X$, $x\in X$ is a fixed
 point, $\varepsilon>0$ is a constant which we use for regularization and eventually let $\varepsilon \to $0; $\xi$ are parameters of hyperplane sections. Let us remark that the 1st factor in the form differs from the form in the definition of $\hat f$ only by the exponent in the denominator.

 \begin {proposition} The form $\kappa_x[f]$ is closed on $X\times {\mathbb R}^n_\xi$.
\end {proposition}

Indeed, in $u$, this a maximal degree form; in $\xi$, we take
\[
a(\xi)=\frac {u+x}{(\langle \xi , (u-x)\rangle -{\rm i} \varepsilon)^{n-1}}
\] and then apply our formula for the differential of the determinant above. The result will contain the factor
\[
\square_{p,q}(u)-\square_{p,q}(x)=1-1=0.
\]

We will reconstruct $f(x)$ by integrating this closed form $\kappa_x[f]$ along different cycles and then regularizing as
 $\varepsilon\to 0 $. Let $\gamma(x)$ be the cycle of sections passing through a fixed point $x\in X$, homological to the sphere $S^{n-2}$. Let $\gamma_0(x)$ be the cycle of sections that pass through $x$ and 0; we call $\gamma_0(x)$ {\em geodesic}.

\begin{proposition}\label{proposition2} For $\varepsilon=0$ we have
\[
 \int_{X\times \gamma(x)} \kappa_x[f]=cf(x),\qquad c=\frac {2(2\pi {\rm i})^{n-1}}{(n-1)!}.
 \]
Here $\varepsilon=0$, $\gamma(x)$ is an $(n-2)$-dimensional cycle in~$\xi$, homologous to the sphere $S^{n-2}$.
\end{proposition}

For the proof it is enough to consider one point $x$ (using the invariance of $X_{p,q}$) and the geodesic cycle $\gamma_0(x)$. Let $x=(1,0,\dots,0)$ and consider $\gamma_0(x)$. On this cycle $\xi_1=0$, $\xi=(0,\eta)$. As a result, the factor that contains the determinant is
\[ (1+u_1)\big[\eta,{\rm d}\eta^{\{n-2\}}\big],
\] and therefore the integral is the sum of two terms.

 If $f(u)$ is even then the 1st term is even, but the 2nd one includes~$u_1f$ and is odd; so the integral of the 2nd one is zero. The integral of the 1st one is $cf(x)$ as follows from the projective Radon inversion formula~\cite{3} in the Cauchy form  (the central projection). If~$f$ is odd, then the 1st term gives zero, but the 2nd one gives $cx_1f(x)=cf(x)$. So the two terms in the form work for even and odd components of the function~$f$ respectively.

For the inversion of $\hat f$ for general cycles $\gamma(x)$ we need to express the restriction of the form~$\kappa_x[f]$ to the cycle through~$\hat f$. It follows from the last proposition that generically we can do it using restrictions of~$\hat f$ and its 1st derivatives (the Cauchy problem for the ultrahyperbolic equation!). In integral geometry we are interested in the special case where it is sufficient to
know restrictions of $\hat f$ -- the characteristic cycles~$\gamma(x)$ (Goursat problem). The description of such cycles (also known as admissible sets) is a very interesting problem, see~\cite{3}, but we investigate here only one important class~-- the horospherical cycles.

\section[The hyperbolic space $X_{1,n-1}$]{The hyperbolic space $\boldsymbol{X_{1,n-1}}$}

Let us start with the case of the two-sheeted hyperboloid $X=X_{1,n-1}$, which carries an action of the group ${\rm SO}(1,n-1)$, and we will only work with one of its sheets: $X_+=\{x_1>0\}$.

The horospherical inversion formula is well known in the hyperbolic case, but we want to illustrate in this example the method that we discuss in this paper. I believe it explains the old observation made by Gelfand that this inversion coincides with Radon's inversion formula.

So we consider the hyperplane sections $L(\xi,p)$ of $X_+$ by the hyperbolic spheres
\[ \langle \xi,x\rangle =p,\qquad  \xi \neq 0.
\]  For $p=0$ (hyperplanes passing through 0) we have hyperbolic geodesic hyperplanes. Hyperbolic geodesic Radon transform is projectively equivalent to the Euclidean Radon transform~\cite{3}.

The sections $L(\xi,p)$ by isotropic hyperplanes
\[ \square _{1,n-1}(\xi)=0, \qquad \xi\neq 0,\]  are called horospheres $E(\xi,p)$. They can be interpreted as limits of hyperbolic spheres. So such~$\xi$ are points of the cone $\Xi_{1,n-1}$, $\xi\neq 0$ (also known as the asymptotic cone) that corresponds to the hyperboloid~$X$. Here $p\neq 0$, since for $p=0$ the sections have no points. In the case of horospheres we call $\hat f$ the horospherical Cauchy transform:
\[ {\mathcal H}f(\xi,p)=\int _X \frac {f(x)}{\langle \xi,x \rangle -p-{\rm i} 0} \omega_{1,n-1}(x,{\rm d}x),\qquad \square_{1,n-1}(\xi)=0.
\] It is sufficient to consider $p=1$. We want to reconstruct $f(x)$ from ${\mathcal H}f$.

Let $\gamma_1(x)$ be the $(n-2)$-dimensional cycle of horospheres passing through $x$. It is enough to consider a fixed point, and so let $x=(1,0,\dots,0)$. Then on $\gamma_1$ we have $\xi_1=p$, $\xi=(p,\eta)$, $p^2=\bigtriangleup(\eta)=(\eta_1)^2+\cdots+(\eta_{n-1})^2$. Therefore
\[ \xi=\big(\sqrt{\bigtriangleup (\eta)},\eta\big),\qquad  p=\sqrt{\bigtriangleup (\eta)}.
\] There is a natural retraction of the horospherical cycle~$\gamma_1(x)$ on the geodesic cycle $\gamma_0(x)$ through cycles $\gamma_\rho(x)$, $0\leq \rho \leq 1$, of the hyperplane intersections $L(\xi,p)$ with
\[
\xi=\big(\rho\sqrt{\bigtriangleup (\eta)},\eta\big), \qquad p=\rho\sqrt{\bigtriangleup (\eta)}.\]
 All $\gamma_\rho$ pass through~$x$. Geometrically, we transform geodesic hyperbolic hyperplanes through hyperbolic spheres in hyperbolic horospheres.

To obtain the horospherical inversion formula we take the cycle $\gamma_1$ in Proposition~\ref{proposition2} and we need to express the restriction of the form~$\kappa_x[f]$ to this cycle in terms of the restriction of the horospherical transform ${\mathcal H}$. We need to remove~$u$ from the determinant $\big[x+u,\xi,{\rm d}\xi ^{\{n-2\}} \big] $ (generically, only this part requires transversal derivatives of~$\hat f$). We can do this specifically for horospheres by using one more trick with determinants.

We want to rewrite formulas below for arbitrary $x\in X$ (not only for $x=(1,0,\dots,0)$). So $\square (x)=1$, $\square (\xi)=0$, $\langle \xi,x\rangle =p$. Here $p \neq 0$. If $\xi=x+\eta$, $\langle \eta,x\rangle=0$, then $\square (\eta)=\langle \eta,\eta\rangle=p-1$.

The transformation of the form $\kappa_x[f]$ uses the following important fact. If $\square (\xi)=0$, $\langle \lambda,\xi\rangle \neq 0$ then the form
\[ \frac {\big[\lambda ,\xi,{\rm d}\xi ^{\{n-2\}}\big]} {\langle\lambda , 	\xi\rangle}\]
is independent of $\lambda$. If we let for simplicity $\xi_1 \neq 0$, then it is enough to add to the 1st row the $j$-th row (any $j$ will do) with coefficient $\xi_j / \xi_1$. As a result, all entries of the 1st row, except for the 1st one, will be 0, while the 1st one will equal
 $\langle\xi,\lambda\rangle/\xi_1$, and thus $\lambda$ disappears from the formula. The conceptual meaning of the $\lambda$-independence is that this form is the residue-form of
 \[ \big[\xi,{\rm d}\xi ^{\{n-1\}}\big]/\square(\xi)\]  on the horospherical cycle
$\square (\xi)=0$.
 A special case of this independence gives
\[ \big[x+u,\xi,{\rm d}\xi^{\{n-2\}}\big]=\frac {\langle\xi, (x+u)\rangle}{\langle\xi ,
	x\rangle}\big[x,\xi,{\rm d}\xi ^{\{n-2\}}\big]=\left(2+\frac {\langle \xi, (u-x)\rangle}{p}\right)\big[x,\xi,{\rm d}\xi ^{\{n-2\}}\big];\]
in the last determinant we can replace $\xi$ with $\eta$.
The substitution in $\kappa_x[f]$ on the horospherical cycle gives
\begin{gather*}
\int _{X \times \gamma_1(x)} f(u) \left(\frac 2 {(\langle \xi , (u-x)\rangle-{\rm i}0)^{n-1}}
+\frac 1 {p(\langle \xi , (u-x)\rangle-{\rm i}0)^{n-2}} \right)\\
\qquad{}\times  \big[u,{\rm d}u^{\{n-1\}}\big] \wedge \big[x,\xi,{\rm d}\xi ^{\{n-2\}}\big]=\frac {2(2\pi {\rm i})^{n-1}}{(n-1)!}f(x).
\end{gather*}
Let us emphasize that in these transformations it was essential that we consider horospheres ($\square (\xi)=0$). Now we can express the form $\kappa$ through the horospherical transform ${\mathcal H}f$. Let us define the differential operator
\[ {\mathcal L}_p=\frac {n-2}p
\left(\frac {\partial}{\partial p}\right)^{n-3} +2\left( \frac {\partial}{\partial p}\right)^{n-2}.
\]
	
\begin{theorem} One has the horospherical inversion formula
\[ f(x)=\frac{n-1}{2(2\pi {\rm i})^{n-1}}\int_{\gamma_1(x)} {\mathcal L}_p\hat f(x+{\rm i}\eta,p)|_{p=1} \big[x,\eta, {\rm d}\eta ^{\{n-2\}}\big].
\]
\end{theorem}

Often in the definition of the hyperbolic horospherical transform, a factor is added which makes the inversion formula slightly different. For the transition from the hyperbolic horosphe\-rical transform to the harmonic analysis on the hyperbolic space $X_{1,n-1}$ -- the spherical Fourier transform on $X$ -- we need to compose the horospherical transform with the Mellin transform on the one-parametric dilations of the cone $\Xi_{1,n-1}$:
\[ \xi \to\lambda \xi,\qquad \lambda \in {\mathbb R}^+.\]

\section[The sphere $S^{n-1}=X_{n,0}$]{The sphere $\boldsymbol{S^{n-1}=X_{n,0}}$}

We use in this section $\Delta$ instead $\square _{n,0}$. Therefore on the $(n-1)$-dimensional sphere $S^{n-1}$ we have $\Delta(x)=1$. We want to find $S^{n-1}$-analogs of the constructions from the last section. However, there are no real horospheres on the sphere since $\Delta(\xi) \neq 0$ if $\xi \neq 0$, and we consider complex horospheres instead~\cite{4}.
In other words, we continue to work with the real sphere and its points, but we take sections by complex hyperplanes $\langle \zeta,z\rangle=p+{\rm i}q$, $\zeta=\xi+{\rm i}\eta$, $\zeta \neq 0 \in {\mathbb C}^n$ with complex points $z\in {\mathbb C}S^{n-1}$, including the complex horospheres $E(\zeta,p+{\rm i}q)$ where $\Delta(\zeta)=0$, $\zeta\neq 0$; the set of such $\zeta$ will be denoted by ${\mathbb C}\Xi$.

We shall focus on complex horospheres (and other hyperplane sections) without real points on $S^{n-1}$. The reason for doing so is that then the Cauchy--Radon transform
\[
\hat f(\zeta,p+{\rm i}q)=\int _S \frac {f(x)}{\langle \zeta,x \rangle -p-{\rm i}q} \omega(x,{\rm d}x)
\]
(where we integrate on the invariant form $\omega_{n,0}$) makes perfect
sense as there are no singularities in the denominator. In this case $\hat f$ is holomorphic in the interior of their domain.

Let us investigate when complex horospheres have no real points on the sphere. For $\zeta=\xi+{\rm i}\eta$ the condition $\zeta \in{\mathbb C}\Xi (\Delta(\zeta)=0)$, $\zeta\neq 0,$ is equivalent to the conditions
\[
\Delta(\xi)=\Delta(\eta), \qquad \langle \xi,\eta\rangle =0.
\] On ${\mathbb C}\Xi$ we have $\Delta(\xi)=\Delta(\eta)>0$. So, on the horospheres, one has $p+{\rm i}q\neq 0$, and we can always take $p+{\rm i}q=1$.

With this normalization the set $\Xi _+$ of horospheres without real points is characterized by the condition
\[ 0<\Delta(\xi)=\Delta(\eta)<1.
\] Indeed, using the action of the orthogonal group, the verification can be reduced to $\xi=(a,0,\dots,0)$, $\eta=(0,b,0,\dots,0)$, in which case there are no real points on the horosphere iff $0\leq |a|^2=|b|^2 \leq 1$. Without the condition $p+{\rm i}q=1$ we have $0<\Delta(\xi)<|p+{\rm i}q|^2$.

On the boundary $\partial \Xi_+$ of the domain $\Xi_+$ we have \[ \Delta(\xi) =\Delta(\eta)=1.\] If $x\in S$ is a (real) point on the horosphere $E(\zeta)$, then $\langle \xi,x\rangle=1$,
$\langle \eta,x\rangle =0$, thus $x=\xi$, and it is the unique real point~$x$ on~$S$; so
\[  \zeta= x+{\rm i} \eta, \qquad \Delta(\eta)=1, \qquad \langle \eta,x\rangle =0.\]
 Therefore $\partial \Xi_+$ projects on the sphere $S^{n-1}$ with fibers the orthogonal $(n-2)$-dimensional spheres.

 Similarly to the construction in the last section, let us consider the geodesic cycle $\gamma_0(x)$ of hyperplanes passing through the center~0 and $x\in S^{n-1}$. We construct a deformation of $\gamma_0(x)$ in the cycle $\gamma_1(x)$ of complex horospheres with the unique real point $x$ as follows. Using the invariance, it is enough to consider $x=(1,0,\dots,0)$. Then the geodesics from $\gamma_0(x)$
are sections of~$S$, by
\[ \big\langle \zeta^0,u\big\rangle =0, \qquad \zeta^0=\big(0,\lambda,\lambda \neq 0 \in {\mathbb R}^{n-2}\big).\]

For $0\leq \delta \leq 1$, let \[
\zeta ^
\delta=\big({\rm i}\delta \sqrt {\Delta(\lambda}),\lambda\big),\qquad p^\delta={\rm i}\delta \sqrt {\Delta(\lambda)},\]
 and $L\big(\zeta^\delta,p^\delta\big)$ are the corresponding sections by the hyperplanes. For $\delta=0$ we have the cycle~$\gamma_0(x)$ of geodesic sections; for $\delta=1$ we have the cycle $\gamma_1(x)$ of complex horospheres, and we have the intermediate cycles $\gamma_\delta(x)$ for $0<\delta<1$.

Sections from all these cycles for $0<\delta \leq 1$ have the unique real point
 $u=x=(1,0,\dots,0)$. Indeed, if $u$ is a real point, then, taking the imaginary parts of both parts of the equation, we have $u_1=1$, but on $S$ there is only the point $u=(1,0, \dots,0)$ with $u_1=1$. So, with the exception of the point $x\in S$, the deformation of the real geodesic planes on~$S$, passing through~$x$, in horospheres tangent to $S$ in $x$, lies in the complex domain of~${\mathbb C}S$.

If we keep the same $\zeta^\delta$ and take $p^\delta(\varepsilon)=p^\delta+{\rm i}\varepsilon$, $\varepsilon>0,$ then the ``shifted'' cycles $\gamma_\delta(x,\varepsilon)$ will consist of sections without real points. Correspondingly, $\gamma_1(x,\varepsilon)$ consists of complex horospheres without real points, and we can consider the horospherical Cauchy transform. On $\gamma_1(x)$ we define it by the regularization for $\varepsilon \to 0$ (boundary values ${\mathcal H}f(x+{\rm i}\eta)$ on $\partial \Xi_+)$. We can now use the fundamental form and repeat all constructions of the last section.
	
\begin{theorem} One has the horospherical inversion formula
\[ f(x)=\frac{n-1}{2(2\pi)^{n-1}}\int_{\gamma_1(x)} {\mathcal L}_p{\mathcal H} f(x+{\rm i}\eta,p)|_{p=1} \big[x,\eta, {\rm d}\eta
^{\{n-2\}}\big].
\]
\end{theorem}

Finally, let us explain how to reproduce the classical version of the harmonic analysis on the sphere in the language of spherical polynomials. The central observation here is that the set~$\Xi_+$, along with its boundary, are invariant relative to the circle action
\[
\zeta \mapsto \exp({\rm i}\theta)\zeta, \qquad \zeta\in \partial\Xi_+.
\] The composition with the horospherical transform gives the decomposition on irreducible components and the inversion gives on the sphere $S$ the decomposition on harmonic polynomials.

\section[The hyperboloid $X_{2,n-2}$]{The hyperboloid $\boldsymbol{X_{2,n-2}}$}

The case $p=2$ contains some important special subcases: for $q=2$ we have the group ${\rm SL}(2;{\mathbb R})$. The imaginary hyperbolic plane ($p=1$, $q=2$) serves as a model of the principal series of representations of this group (all representations appear with multiplicity~1.) For our considerations it is the simplest case where we need to consider both real and complex horospheres. This forces us to combine the technology which we discussed for the hyperbolic space and for the sphere~\cite{2}.

Following our conceptual picture, we are interested in the set $\varTheta$ of horospheres~$E(\zeta,p)$ without real points on $X$, since for them the horospherical Cauchy transform ${\mathcal H}(\zeta,p)$ is well defined. Here $\zeta=\xi+{\rm i}\eta \in {\mathbb C}\Xi$ is nonzero and satisfies $\square\zeta=0$. Using the homogeneity of ${\mathcal H}$, it is enough to consider $p=0,1$. Of principal interest are the boundary horospheres from $\partial \varTheta$: they have intersections with $X$ which disappear under some small perturbations. The horospherical transform ${\mathcal H}$ for them is defined by regularization (boundary values) from $\varTheta$. The set of such horospheres is a dual object to~$X$ and we use them to write the inversion formula.

The first class of horospheres which lie in $\partial \varTheta$ is the class $\varTheta_ R$ of the real horospheres $E(\xi,p)$ for which $\square (\xi)=0$, $\eta =0$, $p\in {\mathbb R}$. Their intersections with $X$ are paraboloids which degenerate to the cone if $p=0$. These horospheres are the boundary for $\varTheta$ since the horospheres $E(\xi,p+{\rm i}\varepsilon)$ have no real points for any $\varepsilon \neq 0$.

For the next class of horospheres $\varTheta_I^\pm$, we suppose that
\[ \square (\xi) \neq 0,\] and this class has two connected components. There are a few degenerate horospheres when $\square (\xi)=\square (\eta)=0$, but we do not need them for our aims. We know that the condition $\square_{1,n-1}(\zeta)=0$ is equivalent to
\[ \square (\xi)=\square (\eta),\qquad \langle \xi,\eta\rangle =0.
\]  With each point $x\in X$ we associate the set of those horospheres $E(\xi+{\rm i}\eta,1)$ that intersect the hyperboloid in the unique point~$x$:
\[ \langle x+{\rm i}\eta,x\rangle  =1, \qquad \langle  \eta,x\rangle  =0.\]

Such horospheres lie on the boundary $\partial \varTheta$, and
to see this let us consider various possibilities. Let $p>0$; then we can take $p=1$. Let \[ \square(\xi)=\square(\eta)>0.\] Then, acting by the group ${\rm SO}(2,q)$, we can make $\xi=(\lambda,0,\dots,0)$, $\lambda >0$, $\eta=(0,\pm\lambda,0,\dots,0)$. If $x\in X$ then $x_1=1/\lambda$, $x_2=0$, and
\[ (x_3)^2+\cdots +(x_n)^2=-1+1/\lambda^2.\] It is impossible to have $\lambda>1$, but for $\lambda=1$ the point $(1,0,\dots,0)$ will be the unique intersection point of the hyperboloid and the horospheres. This exactly means that the horosphere $E(\zeta,1)$ is on the boundary of~$\varTheta$. Apparently, for $\lambda<1$ there is more than one intersection point.

For such $\zeta$ and $p=0$ in the intersection it must be the case that $x_1=x_2=0$, and this is impossible. So such horospheres are interior to $\varTheta$.

If $\square (\xi)<0$ then the canonical form of such $\zeta$ is $(0,0,\lambda,{\rm i}\lambda,0,\dots,0)$, and such horospheres have infinite real intersections for all $p$. The case $\xi=0$, $\eta\neq 0$ can be considered similarly.

The basic consequence of these computations is the following: The horospheres $E(\zeta,p)$ have no real points if
\[ \square (\xi)=\square (\eta)>|p|^2>0;\]
 if $ \square (\xi)=|p|^2>0$, then $x=\xi/p$ has a unique intersection point.
So passing through $x\in X$ are the horospheres $E(x+{\rm i}\eta,1)$ with
\[ \square (\eta)=1,\qquad \langle \eta,x\rangle =0.\] These horospheres are parameterized by the points of the $(n-2)$-dimensional, two-sheeted hyperboloid $\{\square _{1,n-2}(\eta)=1\}$. Now we have enough horospheres for the inversion of the horospherical Cauchy transform. This set $\gamma_1(x) \subset \partial \varTheta$ has three connected components: one component from the real horospheres $\gamma_1^0(x)$ and two components from the complex horospheres~$\gamma_1^\pm(x)$.

Our approach has the same structure as in the examples above. We start from the geodesic hyperplane sections $L(\xi)$ by the hyperplanes $\langle \xi,u\rangle=0$. We fix a point $x\in X$. It is sufficient to take $x=(1,0,\dots,0)$. For these sections the inversion follows from the Radon inversion formula by the projective equivalence. Let us take the cycle $\gamma_0(x)$ of the geodesic sections passing through~$x$. Then $\xi_1=0$ and we have sections
\[ \big\langle \zeta^0,u\big\rangle=0, \qquad \zeta^0=(0,\lambda),\qquad \lambda\neq0\in {\mathbb R}^{n-1}.\]
 For normalization we can take $\lambda$ from the unit sphere~$S^{n-2}$.
We construct the deformation of the geodesic cycle $\gamma_0(x)$ to the horospheric cycle~$\gamma_1(x)$.

For $0 \leq \delta \leq 1$ consider
\[ \zeta ^\delta=\big({\rm i}\delta \sqrt {\square (\lambda)},\lambda\big), \qquad p^\delta={\rm i}\delta \sqrt {\square(\lambda)}\] and the cycles $\gamma_\delta(x)$ of the hyperplane sections $L\big(\zeta^\delta,p^\delta\big)$ by the hyperplanes $\big\{\big\langle \zeta^\delta,u\big\rangle=p^\delta \delta\big\}$ that depend on $\lambda \in S^{n-2}$. For $\delta =0$ we have the geodesic cycle~$\gamma_0(x)$; for $\delta=1$ we have the horospherical cycle $\gamma_1(x)$ that consists of the real horospheres for $\square_{1,n-1}(\lambda)<0$ and the 2-component set of complex horospheres. Let $\gamma_1^0(x)$, $\gamma^\pm_1(x)$ be these three connected components of the horospherical cycle. Their joint boundary consists of degenerated horospheres. Thus we have constructed the deformation of the geodesic cycle in the horospheric cycle through the intermediate cycles $\gamma^\delta (x)$, $0<\delta<1$ of hyperplane sections (not horospheres).

Correspondingly, in $\varTheta$ we have the two-component domain $\varTheta^\pm$ of horospheres $E(\zeta,1)$ with
\[\square (\xi)=\square (\eta)>1\] on whose boundaries lie $\gamma^\pm(x)$. The horospherical transform is holomorphic in the domains $\varTheta^\pm$. Therefore the horospherical transform has 3 components: ${\mathcal H}_ R(\xi)$, ${\mathcal H}_I^\pm(\zeta)$. We can construct the inversion of the horospherical transform on $\gamma_1(x)$ using the fundamental form. We will discuss the details at another time.

The spherical Fourier transform also has 3 components corresponding to 3 components of the horospherical transform. The first one is connected with the action of multiplicative real group~${\mathbb R}^\times$ on the set~$\varTheta_R$ of real horospheres
\[\xi\mapsto \lambda\xi,\] and, as in the hyperbolic case, we take the composition of the horospherical component~${\mathcal H}_R$ and the classical Mellin transform. It gives the projection on the continuous principal series.

The domains $\varTheta^\pm _I$ are invariant relative to the complex semigroup: $\zeta \mapsto \lambda \zeta$, $\lambda \in {\mathbb C}$, $|\lambda| >1$. The decomposition of ${\mathcal H}^\pm$ in the Taylor series in $\lambda$ gives projections on the holomorphic and antiholomorphic discrete series of representations.

\pdfbookmark[1]{References}{ref}
 \LastPageEnding

\end{document}